\input amstex
\documentstyle{amsppt}
\magnification=1200
\UseAMSsymbols
\NoBlackBoxes
\topmatter
\title On Complemented Subspaces of Sums and Products of Banach spaces
\endtitle
\rightheadtext{COMPLEMENTED SUBSPACES}
\author M.I.Ostrovskii
\endauthor
\address Mathematical Division, Institute for Low Temperature Physics
and Engineering, 47 Lenin avenue, 310164 Kharkov, UKRAINE \endaddress
\email mostrovskii\@ilt.kharkov.ua \endemail
\abstract
It is proved that there exist complemented subspaces
of countable topological products (locally convex direct sums) 
of Banach spaces which cannot be
represented as topological products (locally convex direct sums) 
of Banach spaces
\endabstract
%\thanks   \endthanks
\endtopmatter

\document
%\UseAMSsymbols
The problem of description of complemented subspaces of a given locally
convex
space is one of the general problems of structure theory of locally 
convex spaces.
In investigations of this problem in the particular case of
spaces represented as countable products of Banach spaces (see \cite{D1},
\cite{D2},
\cite{DO}, \cite{MM1}) the following problem arose (see
\cite{D2, p.~71}, \cite{MM2, p.~147}):
Is every complemented subspace of a topological 
product (locally convex direct sum) of a countable family
of Banach spaces isomorphic to a topological product (locally 
convex direct sum) of Banach spaces?
G.Metafune and V.B.Moscatelli \cite{MM3, p.~251} conjectured that this is
false,
in general. The purpose of the present note is to prove this conjecture.

Our sources for basic concepts and results of Banach space theory
and the theory of topological vector spaces are, 
respectively, \cite{LT} and \cite{RR}.

Let us fix some terminology and notation.
The algebra of all continuous linear operators on a Banach space $X$ will
be denoted by $L(X)$.
The identity mapping of a linear space $W$ is denoted by $I_W$. Let
$\{X_n\}_{n=1}^\infty$ be a sequence of Banach spaces. We denote
their Cartesian product endowed with the product topology by
$\prod_{n=1}^\infty X_n$ and call it {\it topological product}.
We denote the locally convex direct sum of spaces $\{X_n\}_{n=1}^\infty$ 
by $\oplus_{n=1}^\infty X_n$. A linear subspace $Y$ of a topological
vector space $Z$ will be called {\it complemented} if there is a
continuous linear mapping $P$ of $Z$ onto $Y$ such that $P^2=P$.
If $B$ is a subset of linear space $V$, then the linear subspace
of $V$ generated by $B$ will be denoted by lin$B$. The dual of
a locally convex space $Z$ endowed with its strong topology
will be denoted by $Z'_\beta$.
\proclaim{Theorem 1} A. There exists a sequence $\{X_n\}_{n=1}^\infty$
of Banach spaces and a complemented subspace $Y$ in $X=\oplus_{n=1}^\infty
X_n$, such that $Y$ is not isomorphic to a locally convex
direct sum of Banach spaces.

B. There exists a sequence $\{Z_n\}_{n=1}^\infty$
of Banach spaces and a complemented subspace $W$ in $Z=\prod_{n=1}^\infty
Z_n$, such that $W$ is not isomorphic to a topological 
product of Banach spaces.
\endproclaim
\demo{Proof} First we shall prove part A.
The main tool of our construction is the space with the property of
bounded approximation but without $\pi$-property, constructed by C.J.Read
\cite{R}. Now we describe those
details of Read's construction which we will use.

Let $V_0$ be a vector space of countable dimension with basis
$\{e_i\}_{i=1}^\infty$. By ${\Bbb N_6}$ we denote the set
${\Bbb N}\backslash\{1, 2, 3, 4, 5\}$. Let $a_r=5\cdot10^{r-2},\
b_r=10^r,\ r\in{\Bbb N_6}$ and $V=$lin$\{e_i:\ i>a_6\}$.
Let us introduce in $V$ three collections of finite-dimensional subspaces:
$$V_r=\hbox {lin}\{e_i:\ a_r<i\le b_r\},\ r\in{\Bbb N_6},$$
$$W_r=\hbox {lin}\{e_i:\ a_{r+1}<i\le b_r\}=V_r\cap V_{r+1},\ r\in{\Bbb N_6},$$
$$U_r=\cases \hbox {lin}\{e_i:\ b_{r-1}<i\le b_r\},\
r\in{\Bbb N_6}\backslash\{6\},\\ V_6,\ r=6.
\endcases
$$
We introduce the notation
$$E(r)=\cases \{i:\ b_{r-1}<i\le b_r\},\
r\in{\Bbb N_6}\backslash\{6\},\\\{i:\ a_6<i\le b_6\},\ r=6.
\endcases
$$

By $\pi_r^V:V\to V_r,\ \pi_r^W:V\to W_r,\ \pi_r^U:V\to U_r$ we denote
natural projections. For example,
$$\pi^V_r(\sum_{i>a_6}d_ie_i)=\sum_{a_r<i\le b_r}d_ie_i.$$

Given some norms $\{||\cdot||_{(r)}\}_{r=6}^\infty$ on
$V_r$, we endow $V$ with the norm
$$||x||^2=\sum_{r=6}^\infty[||\pi_r^V(x)||^2_{(r)}+
\frac1{\log r}(||\pi_r^W(x)||^2_{(r)}+||\pi^W_r(x)||^2_{(r+1)})]$$

Let $R$ be the completion of $V$ under the norm $||\cdot||$.

Let $A$ be a subset of ${\Bbb N_6}$. We denote
the closure in the norm topology of the subspace of $R$ spanned by vectors
$$\{e_i:\ i\in\bigcup_{j\in A}E(j)\}$$
by $R_A$.
The arguments of C.J.Read \cite{R} imply the
following result.
\proclaim{Theorem 2} For some collection of norms $||\cdot||_{(r)},\
(r\in{\Bbb N_6})$ there exist a convergent to zero sequence of
positive numbers $\{\varepsilon(r)\}_{r=6}^\infty $ and numbers
$0<\alpha,\beta<\infty$ such that for every $A\subset {\Bbb N}_6$
and every $j\in A$ there exists a mapping $\xi_j:L(R_A)\to
\Bbb R$ such that the following conditions are satisfied:
\roster
\item"(a)" If $\ \pi^W_jT|_{W_j}=\mu I_{W_j}$, then $\xi_j(T)=\mu$.

\item"(b)" $|\xi_j(TS)-\xi_j(T)\xi_j(S)|\le\alpha j^{-\beta}||T||||S||.$

\item"(c)" If $j,j+1\in A$, then $|\xi_{j+1}(T)-\xi_j(T)|
\le \varepsilon(j+1)||T||.$
\endroster
\endproclaim

Starting from this point we suppose that $R$ is the space constructed
by the described method using some collection of norms
$\{||\cdot||_{(r)}\}_{r=6}^\infty$ satisfying the conditions of
Theorem 2. Precise description of norms $||\cdot||_{(r)}$
does not matter for us. We shall use only the definition of
$||\cdot||$ and the properties of spaces $R_A$ listed in
Theorem 2.

Let $f:\Bbb N_6\to\Bbb R$ be some function.
Each vector $v\in V$ can be in a unique manner represented as follows
$$v=\sum_{j=6}^\infty\sum_{i\in E(j)}d_ie_i.$$
Therefore the formula
$$M_fv=\sum_{j=6}^\infty\sum_{i\in E(j)}f(j)d_ie_i$$
defines a linear mapping on $V$.
\proclaim{Lemma 1} If function $f$ is such that for some
$C<\infty$ and each $k\in \Bbb N_6$ the conditions
$$|f(k)|\le C,$$
$$|f(k)-f(k+1)|\le C/\sqrt{\log k},$$
are satisfied, then $M_f$ is bounded with respect to the norm
$||\cdot||$, and $||M_f||\le2C.$
\endproclaim

\demo{Proof}Let us note that for every $v\in V$ we have
$$\pi_6^VM_fv=f(6)\pi^V_6v,$$
$$\pi_r^VM_fv=f(r)\pi_r^Vv+(f(r-1)-f(r))\pi^W_{r-1}v,\
r\in\Bbb N_6\backslash\{6\}.$$
$$\pi^W_rM_fv=f(r)\pi^W_rv,\ r\in\Bbb N_6.$$
Hence, we have
$$||M_fv||^2=$$
$$\sum_{r=6}^\infty[||\pi_r^VM_fv||^2_{(r)}+
\frac1{\log r}(||\pi_r^WM_fv||^2_{(r)}+||\pi^W_rM_fv||^2_{(r+1)})]\le$$
$$|f(6)|^2||\pi_6^Vv||^2_{(6)}+$$
$$\sum_{r=7}^\infty(|f(r)|||\pi_r^Vv||_{(r)}+|f(r)-f(r-1)|
||\pi^W_{r-1}v||_{(r)})^2+$$
$$\sum_{r=6}^\infty |f(r)|^2\frac1{\log r}(||\pi^W_rv||^2_{(r)}+
||\pi^W_rv||^2_{(r+1)})\le$$
$$C^2||\pi^V_6v||^2_{(6)}+\sum_{r=7}^\infty2C^2||\pi^V_rv||^2_{(r)}+
\sum_{r=7}^\infty\frac{2C^2}{\log(r-1)}||\pi^W_{r-1}v||^2_{(r)}+$$
$$\sum_{r=6}^\infty\frac{C^2}{\log r}(||\pi^W_rv||_{(r)}^2+
||\pi_r^Wv||_{(r+1)}^2)\le$$
$$4C^2||v||^2.$$
The required inequality follows.
\enddemo

Therefore, if the conditions of Lemma 1 are satisfied then $M_f$
may be considered as an element of $L(R)$.

Let $\{f_n\}_{n=1}^\infty$ be a sequence of functions
satisfying the conditions of Lemma 1.
Introduce notation $M_n:=M_{f_n};\ A(n)$:=supp$ f_n\subset\Bbb N_6$.
\proclaim{Lemma 2} There exists a sequence $\{f_n\}_{n=1}^\infty$ 
such that
\roster
\item"(a)" $\sup_n||M_n||<\infty.$

\item"(b)" If $j\in A(n)$, then $f_{n+1}(j)=1$.

\item"(c)" For every $j\in\Bbb N_6$ there exists $m\in\Bbb N$ such that
$f_m(j)=1.$

\item"(d)" For every $n\in\Bbb N$ and every
$r\in\Bbb N_6$ there exists $j>r$ such that $j\in A(1)$ and
there exists $i>j$ such that $i\notin A(n)$ and
$\{j,j+1,\dots,i\}\subset A(n+1)$.
\endroster
\endproclaim
\demo{Proof} Let $\{p(s)\}_{s=1}^\infty\subset\Bbb N_6$ be an increasing
sequence of natural numbers such that
$$\sum_{k=p(s)}^{p(s+1)-1}\sqrt{\log k}\ge 1.$$
Let $\{s(m)\}_{m=1}^\infty$ be the sequence defined by the equality
$$s(m)=1+2+3+...+m=m(m+1)/2.$$
Let $A(n)\subset\Bbb N_6$ be the sets represented as unions of
intervals of integers in the following way.
$$A(n)=\Bbb N_6\cap(\cup_{m=1}^\infty\{p(s(m)-n)+1,\dots,p(s(m)+n)-1\}).$$

From the definitions of the sequences $\{p(s)\}_{s=1}^\infty$ and
$\{s(m)\}_{m=1}^\infty$ it follows that there exists a sequence
$\{f_n\}_{n=1}^\infty$ of functions,
$f_n:\Bbb N_6\to\Bbb R$, satisfying the following conditions.
\roster
\item"(e)" $(\forall k\in\Bbb N_6)(0\le f_n(k)\le 1);$
\item"(f)" supp$f_n=A(n);$
\item"(g)" $\{k:\ f_1(k)=1\}=\{p(s(m))\}_{m=1}^\infty;$
\item"(h)" $\{k:\ f_n(k)=1\}=A(n-1)\ (n\ge2);$
\item"(i)" $|f_n(k)-f_n(k+1)|\le\sqrt{\log k}.$
\endroster

By Lemma 1 the conditions (e) and (i) imply (a). It is clear that the conditions
(f)--(h) imply (b) and (c).

In order to verify the condition (d) we let $j$ be any element of the
sequence 
$$\{p(s(m))\}_{m=1}^\infty,$$
$j=p(s(k))$, for which $j>r$ and
$k>2n$. Let $i=p(s(k)+n)$. By the conditions (f) and (g) it follows
that $i$ and $j$ satisfy the condition (d). Lemma 2 is proved.
\enddemo
 
It is clear that the conditions (a), (b) and (c) of Lemma 2 imply that 
the sequence
$\{M_n\}_{n=1}^\infty$ converges to the identity operator in the 
strong operator topology.

Let $C(n)=A(n)\backslash A(n-1)$ (we set $A(0)=\emptyset$). Set 
$$D_n=R_{C(n)}\ (n\in\Bbb N),$$
$$G=R\oplus_p(\sum^\infty_{n=1}\oplus D_n)_p,$$
where $1<p<\infty$ is arbitrary number. 

We recall that by the definitions of direct sums, $G$ is the space
of those sequences of the form
$$g=(x,y_1,y_2,\dots,y_n,\dots),$$
for which $x\in R$, $y_n\in D_n$ and
$$||x||^p+\sum_{n=1}^\infty||y_n||^p<\infty.$$
Linear operations on $G$ are defined in the coordinatewise manner
and the norm is defined by the formula:
$$||g||=(||x||^p+\sum_{n=1}^\infty||y_n||^p)^{1/p}.$$

We introduce operators $S_n\in L(G)\ (n\in\Bbb N)$
by the formulas:
$$S_n(x,y_1,y_2,\dots,y_n,\dots)=$$
$$(M_nx+y_n,y_1,\dots,y_{n-1},(M_n-M_n^2)x+(I-M_n)y_n,0,0,\dots).$$

These operators are well-defined because by the condition (b) of Lemma 2
the image of $(M_n-M_n^2)$ is contained in $D_n=R_{C(n)}$. 

By direct verification it follows that the sequence $\{S_n\}_{n=1}^\infty$
is uniformly bounded, $S_n$ are projections and
$$S_nS_m=S_{min(m,n)}\ (m,n\in\Bbb N).$$
Furthermore,
$$(\forall g\in G)(\lim_{n\to\infty}S_ng=g),$$
where the limit is taken in the norm topology. Therefore
$\{S_n\}_{n=1}^\infty$ generates a Schauder decomposition of $G$.

The construction above is taken from \cite{CK}, it goes back to W.B.Johnson
\cite{J}.

The sequence $\{S_nG\}_{n=1}^\infty$ is an increasing sequence of
subspaces of $G$. Let $X$ be a strict inductive limit of this
sequence. Since each $S_{n-1}G$ is complemented in $S_nG$, then
$X$ is isomorphic to a locally convex direct sum of a sequence
of Banach spaces. Since $X=\cup_{n=1}^\infty S_nG$ then $X$ may
be considered as a subset and even a linear subspace of $G$.

Let us denote by $P$ the projection on $G$ defined by the formula
$$P(x,y_1,y_2,\dots,y_n,\dots)=(x,0,0,\dots,0,\dots).$$

By definition of operator $S_n$ it follows that
$$S_nG\subset R_{A(n)}\oplus_p(\sum_{k=1}^n\oplus R_{C(k)})_p
\oplus_p(\sum_{k=n+1}^\infty\oplus\{0\})_p.
\eqno{(1)}$$
By the condition (b) of Lemma 2 the restriction of $M_n\ (n\ge 2)$
to $R_{A(n-1)}$ coincides with the identity operator. Therefore
$$S_nG\supset R_{A(n-1)}\oplus_p(\sum_{k=1}^{n-1}\oplus R_{C(k)})_p
\oplus_p(\sum_{k=n}^\infty\oplus\{0\})_p.
\eqno{(2)}$$

Therefore $PS_nG\subset S_{n+1}G$ and $P$ may be considered also
as a projection on $X$. Let $Y=PX$. We claim that
\roster
\item"(I)"  $P$ is a continuous projection on $X$.
\item"(II)"  $Y$ considered as a subspace of $X$ is not isomorphic to
a locally convex direct sum of Banach spaces.
\endroster

The topologies
of $X$ and $G$ coincide on each of $S_nG$ \cite{RR, p.127}.
Since $X$ is a strict inductive limit of $\{S_nG\}_{n=1}^\infty$,
then in order to prove statement (I), it is sufficient to prove 
that the restriction of $P$ to each
of $S_nG$ is continuous \cite{RR, p.79}. But this immediately
follows from the inclusion $PS_nG\subset S_{n+1}G$ and the
continuity of $P$ considered as a mapping from $G$ to $G$
(with its initial Banach topology).

Let us prove (II). We have
$$Y=\cup_{n=1}^\infty PS_nG.$$
Therefore
$$R_{A(n-1)}\oplus_p(\sum_{k=1}^\infty\oplus\{0\})_p\subset
Y\cap S_nG\subset R_{A(n)}\oplus_p(\sum_{k=1}^\infty\oplus\{0\})_p.$$
Hence 
$$Y=\cup_{n=1}^\infty (R_{A(n)}\oplus_p
(\sum_{k=1}^\infty\oplus\{0\})_p).$$

Let us show that the topology on $Y$ induced by the strict inductive
topology of $X$ coincides with topology of the strict inductive
limit of the sequence
$$\{R_{A(n)}\oplus_p
(\sum_{k=1}^\infty\oplus\{0\})_p\}_{n=1}^\infty.
\eqno{(3)}$$
We denote
$$R_{A(n)}\oplus_p
(\sum_{k=1}^\infty\oplus\{0\})_p$$
by $K_n\ (n\in\Bbb N)$.

From (2) it follows that the intersection of arbitrary neighbourhood
of zero in $X$ with $K_{n-1}$ contains a neighbourhood of zero in
$K_{n-1}$. Hence the topology induced on $Y$ from $X$ is not stronger
than the strict inductive topology of the sequence (3).

On the other hand, let $\tau$ be a convex neighbourhood of zero in the
strict inductive limit of the sequence (3). Then $\tau_n:=\tau\cap K_n$
is a neighbourhood of zero in $K_n$. Let
$$\theta_n:=\{(x,y_1,\dots,y_n,\dots)\in G:\
(x,0,\dots,0,\dots)\in\tau_n\}.$$
We have $\theta_n\cap Y=\tau_n$. From (1) it follows that
$\sigma_n:=\theta_n\cap S_nG$ is a neighbourhood of zero in $S_nG$. 
It is also clear that $\{\sigma_n\}_{n=1}^\infty$ is an increasing
sequence of convex sets. Therefore $\sigma:=\cup_{n=1}^\infty\sigma_n$
is a neighbourhood of zero in $X$.
It is easy to see that $\sigma\cap Y\subset\tau$. Hence the topology of the
strict inductive limit of the sequence (3) is not stronger than
the topology induced on $Y$ from $X$. Hence these topologies
coincide.

So it remains to prove that the strict inductive
limit of the sequence $\{R_{A(n)}\}_{n=1}^\infty$ is not 
isomorphic to a locally convex direct sum of Banach spaces.

Assume the contrary. Let us denote the strict inductive limit
of $\{R_{A(n)}\}_{n=1}^\infty$ by $S$ and let
$S=\oplus_{\lambda\in\Lambda}S_\lambda$, where $S_\lambda$ are
Banach spaces. Since the space $R_{A(1)}$ is a Banach one, it
has a bounded neighbourhood of zero. By description of bounded
sets in locally convex direct sums \cite{RR, p.92} it
follows that for some finite subset $\Delta\subset\Lambda$
the space $R_{A(1)}$ is contained in
$$T=(\oplus_{\lambda\in\Delta}S_\lambda)\oplus
(\oplus_{\lambda\in\Lambda\backslash\Delta}\{0\}).$$
The space $T$ is isomorphic to a Banach space. Hence it has
a bounded neighbourhood of zero. By description of bounded
sets in strict inductive limits \cite{RR, p.129} it follows
that $T$ is contained in $R_{A(k)}$ for some $k\in\Bbb N$.

It is clear that $T$ is complemented in $S$. Hence $T$
is complemented in $R_{A(j)}$ for every $j\ge k$.
Let $Q\in L(R_{A(k+1)})$ be some projection onto $T$.
Let $r\in\Bbb N_6$ be such that for every $s\ge r$
we have $\varepsilon(s)||Q||<1/10$ and
$\alpha s^{-\beta}||Q||^2<1/100.$ Let $i$ and $j$ be natural
numbers satisfying the condition (d) of Lemma 2 for given
$r$ and for $n=k$. Since $R_{A(1)}$ is in the image of
operator $Q$ and $j\in A(1)$, then by condition (a) of
Theorem 2 it follows that $\xi_j(Q)=1$. Since the
image of $Q$ is in $R_{A(k)}$ and $i\notin A(k)$, then
$\pi_i^WQ=0$. So by the condition (a) of Theorem 2
it follows that $\xi_i(Q)=0$.

Since we have $\{j,j+1,\dots,i\}\subset A(k+1)$ (condition (d)
of Lemma 2), then by the parts (b) and (c) of Theorem 2 we have
$$(\forall q\in\{j,\dots,i-1\})(|\xi_q(Q)-\xi_{q+1}(Q)|<1/10)$$
and
$$(\forall q\in\{j,\dots,i\})(|\xi_q(Q)-\xi_q^2(Q)|<1/100).$$
The last assertion implies that either $|\xi_q(Q)|<1/10$ or
$|\xi_q(Q)|>9/10$. We arrive at a contradiction. Hence $S$
is not isomorphic to a locally convex direct sum of Banach spaces.
The part A of Theorem 1 is proved.

Let us turn to the part B.

Let
$$H=(\sum_{r=6}^\infty\oplus(V_r,||\cdot||_{(r)}))_2
\oplus_2(\sum_{r=6}^\infty\oplus(W_r,||\cdot||_{(r)}))_2
\oplus_2(\sum_{r=6}^\infty\oplus(W_r,||\cdot||_{(r+1)}))_2.$$
From the description of the Banach space $R$ it is clear that
the mapping $T:R\to H$ defined as
$$Tx=(\{\pi^V_r(x)\}_{r=6}^\infty,
\frac1{\sqrt{\log r}}\{\pi^W_r(x)\}_{r=6}^\infty,
\frac1{\sqrt{\log r}}\{\pi^W_r(x)\}_{r=6}^\infty)$$
is an isometric embedding. Because the spaces $V_r$ and
$W_r$ are finite dimensional, then by the well-known rescription
of the duals of direct sums it follows that $H$ is reflexive.
Hence $R$ is reflexive. By the construction of $G$
it follows that $G$ is also a reflexive Banach space.
Hence by well-known description of duals of locally
convex direct sums and topological products \cite{RR, p.93}
it follows that the space $X$ is also reflexive.
Let $Z=X'_\beta$. From the description of a dual of $X$
\cite{RR, p.93} and from the description of bounded
sets of $X$ \cite{RR, p.92} it follows that $Z$ is
isomorphic to a countable topological product of Banach spaces.

Let $P'$ be the adjoint mapping of $P$. It is well-known
\cite{RR, p.48} that $P'$ is continuous in $Z$. Hence $Z$
is isomorphic to a direct sum
$$W\oplus U,$$
where $W=P'Z$ and $U=\{z\in Z:\ P'z=0\}$ (see \cite{RR, p.96}).
Hence
$$Z'_\beta=W'_\beta\oplus U'_\beta,$$
where $W'_\beta$ is the image of $P'':Z'_\beta\to Z'_\beta$.
But because of reflexivity $Z'_\beta=X$ and $P''=P$. 
Hence $Y=W'_\beta$.

So if we suppose that $W$ is isomorphic to a topological
product of Banach spaces, then by description of duals
of topological products \cite{RR, p.93} it would follow
that $Y$ is isomorphic to a locally convex direct sum
of Banach spaces.

The part B of Theorem 1 is also proved.
\enddemo

\enddocument